\documentclass{amsart}
\usepackage[latin1]{inputenc}
\usepackage[english]{babel}
\usepackage{csquotes}

\usepackage{amsmath,amsthm}

\theoremstyle{plain}                    
\newtheorem{theorem}{Theorem}[section]
\newtheorem{lemma}[theorem]{Lemma}
\newtheorem{proposition}[theorem]{Proposition}
\newtheorem{corollary}[theorem]{Corollary}
\theoremstyle{definition}
\newtheorem{definition}[theorem]{Definition}
\newtheorem{example}[theorem]{Example}
\newtheorem{remark}[theorem]{Remark}

\numberwithin{equation}{section}

\usepackage{amssymb}
\usepackage{mathabx}

\newcommand{\ff}{\mathbb{F}}

\newcommand{\rr}{\mathbb{R}}
\newcommand{\qq}{\mathbb{Q}}
\newcommand{\zz}{\mathbb{Z}}
\newcommand{\nn}{\mathbb{N}}

\newcommand{\raag}[1]{A_{#1}} 
\newcommand{\block}[1]{\mathcal B({#1})} 
\newcommand{\vv}[1]{V(#1)} 
\newcommand{\ee}[1]{E(#1)} 

\newcommand{\bldeg}{\operatorname{bldeg}}
\newcommand{\ind}{\operatorname{I}}
\newcommand{\rk}{\operatorname{rk}}

\newcommand{\chars}[1]{S(#1)} 

\usepackage{color} 
\usepackage[colorlinks=true,linktoc=all,linkcolor=blue,citecolor=red,filecolor=pink,urlcolor=blue]{hyperref}
\usepackage[hyperref=true,style=alphabetic,citestyle=alphabetic,backend=bibtex,maxnames=100,doi=false,isbn=false,url=false,giveninits=true]{biblatex}
\bibliography{biblio}

\usepackage{tikz}
\usepackage{float}
\usetikzlibrary{arrows}

\usepackage{todonotes}

\title[Graphical splittings of Artin kernels]{Graphical splittings of Artin kernels}

\author{Danielle  Barquinero}
\email{danielle-barquinero@uiowa.edu}

\author{Lorenzo Ruffoni}
\address{Department of Mathematics - Florida State University, 1017 Academic Way, Tallahassee FL 32306, USA}
\email{lorenzo.ruffoni2@gmail.com}

\author{Kaidi Ye}
\email{deloresye@gmail.com}

\subjclass[2010]{20F36, 20F65, 20E08}
\keywords{Right-angled Artin groups, Bestvina-Brady groups, Artin kernels, Bass-Serre theory, chordal graphs, block graphs, Bieri-Neumann-Strebel invariant.}
 
\begin{document}

\maketitle

\begin{abstract}
We study Artin kernels, i.e. kernels of discrete characters of right-angled Artin groups, and we show that they decompose as graphs of groups in a way that can be explicitly computed from the underlying graph.
When the underlying graph is chordal we show that every such subgroup either surjects to an infinitely generated free group or is a generalized Baumslag-Solitar group of variable rank.
In particular for block graphs (e.g. trees), we obtain an explicit rank formula, and discuss some features of the space of fibrations of the associated right-angled Artin group.
\end{abstract}

\tableofcontents


\section{Introduction}

Given a finite simplicial graph $\Gamma$ the associated right-angled Artin group (RAAG) is the group $\raag \Gamma$ generated by the vertices of $\Gamma$, with two generators commuting if and only if the corresponding vertices are connected by an edge in $\Gamma$; in other words, this is the group with the presentation:
$$\raag{\Gamma} = \langle x \in \vv{\Gamma} \ | \ [x,y] \in \ee{\Gamma} \rangle$$

As such, RAAGs interpolate between free groups and free abelian groups. Moreover they enjoy many good algebraic properties (e.g. they are torsion free, linear over the integers, residually finite, \dots).
Despite their very simple definition, RAAGs turn out to display a surprisingly rich geometric behavior and a very diverse collection of subgroups, which has made them one of the main characters in recent breakthroughs in low-dimensional geometry, as in the work of Haglund and Wise on special cube complexes (see \cite{HW08}) and of Agol on the Virtual Haken Conjecture (see \cite{A13}).\par

One of the most appealing aspects of the theory of RAAGs is that many of their properties can be detected and computed directly from the structure of the underlying graph. For instance a RAAG $\raag \Gamma$ splits as a free product precisely when $\Gamma$ is disconnected, and Clay has shown in \cite{C14} that all non-trivial splittings of $\raag \Gamma$  over $\zz$ correspond to cut vertices of $\Gamma$, and more generally Groves and Hull have shown in \cite{GH17} that $\raag \Gamma$ splits over an abelian subgroup if and only $\Gamma$ is disconnected, or complete, or contains a separating clique (i.e. a full complete subgraphs whose removal disconnects $\Gamma$).\par

This paper deals with fibrations of RAAGs, i.e. non-trivial  homomorphisms 
$$1\to K \to \raag \Gamma \to \zz \to 1$$
with finitely generated kernel. When $\Gamma$ is a tree, $\raag \Gamma$ is known to be the fundamental group of a compact 3-manifold   (see \cite{D87}), and by classical result  of Stallings (see \cite{S62}) fibrations of $\raag \Gamma$ correspond to fibrations of the manifold over the circle, in a way which is well understood in terms of the   Thurston norm on homology (see \cite{T86}). For general graphs these fibrations can be conveniently organized by the so-called Bieri-Neumann-Strebel invariant of $\raag \Gamma$, a tool introduced in \cite{BNS87} as a generalization of Thurston's picture for general finitely generated groups, and explicitly computed for RAAGs by Meier and VanWyk in \cite{MV95}.\par

In \cite{BB97} Bestvina and Brady considered the natural fibration obtained by sending all   the vertices of a connected graph $\Gamma$ to $1 \in \zz$, and studied the corresponding kernel, today known as the Bestvina-Brady group associated to $\Gamma$. They showed that its finiteness properties are completely determined by the topology of $\Gamma$; for instance, the Bestvina-Brady group is finitely presented if and only if the flag complex generated by $\Gamma$ is simply connected. They used this correspondence to show that certain Bestvina-Brady groups provide  counterexamples  either to the Eilenberg-Ganea Conjecture or to the Whitehead Conjecture.  \par

In this paper we consider the generalization of the classical Bestvina-Brady groups given by Artin kernels, i.e. kernels of general homomorphisms $f:\raag \Gamma \to \zz$ (see \cite{BG99,PS09}), and we are interested in exploring their finiteness behavior in terms of their splittings induced by decompositions of the underlying graph $\Gamma$.
A splitting of $\Gamma$ is a decomposition of $\Gamma$ into two full subgraphs $\Gamma=\Gamma_1\cup \Gamma_2$ meeting along a full subgraph $\Gamma_3=\Gamma_1\cap \Gamma_2$; we denote the splitting by $(\Gamma_1,\Gamma_2,\Gamma_3)$ and emphasize that none of these graphs has to be connected. 
Our results concern the ways in which such a splitting can induce a splitting of the Artin kernels associated to $\Gamma$, and are inspired by the work of Cashen and Levitt on free-by-cyclic groups (see \cite{CL16}), in which analogous tame/wild dichotomies and rank formulas are obtained.
Our main result is the following (see Theorem \ref{thm:generaldichotomy} below for more details). 

\begin{theorem}
Let $(\Gamma_1,\Gamma_2,\Gamma_3)$ be a splitting of a connected graph $\Gamma$, let $f:\raag \Gamma \to \zz$ be non-trivial, and let $f_k$ denote the restriction of $f$ to $\raag{\Gamma_k}$ for $k=1,2,3$. Then
\begin{enumerate}
\item  (wild)  if  $f_3= 0$, then $\ker (f)$ surjects onto $\ff_\infty$;
\item (tame) if  $f_3\neq 0$, then $\ker (f)$ splits as a finite graph of groups $\mathcal G$, with $[f(\raag \Gamma) : f(\raag{\Gamma_k})]$ vertex groups isomorphic to $\ker(f_k)$ for $k=1,2$, and $[f(\raag \Gamma) : f(\raag{\Gamma_3})]$ edge groups isomorphic to $\ker(f_3)$.
\end{enumerate}
\end{theorem}

Notice that the two cases here are not 
disjoint from each other, i.e. the same character $f$ can display tame behavior with respect to one splitting, and wild behavior with respect to another (see Remark \ref{rmk:general_dic_caveat}). In order to obtain a complete dichotomy, we restrict to chordal graphs in \S \ref{sec:chordal}, where we prove the following (see Theorem \ref{thm:chordal}).

\begin{theorem} 
Let $\Gamma$ be a connected chordal graph, and let $f:\raag \Gamma \to \zz$ be non-trivial. Then exactly one of the following holds.
\begin{enumerate}
    \item (wild) $\ker(f)$ surjects onto $\ff_\infty$;
    
    \item (tame) $\ker(f)$ is a vGBS group (generalized Baumslag-Solitar group of variable rank), i.e. splits as a finite graph of  finitely generated free abelian   groups. 
\end{enumerate}
\end{theorem}

In the tame case the graph of groups is obtained by applying the previous Theorem to a sequence of splittings that completely decompose $\Gamma$ into the collection of its maximal cliques. It is worth mentioning that the vGBS groups appearing in the tame case are groups for which the multiple conjugacy problem is well understood (see Beeker \cite{B15}).  While it is not known in general if this is the only way an Artin kernel decomposes over abelian subgroups, Chang has recently shown in \cite{C20} that the Bestvina-Brady group over a general graph $\Gamma$ splits over an abelian subgroup if and only if $\Gamma$ is disconnected, or complete, or has a separating clique (the same statement holds for the RAAG $\raag \Gamma $ itself by \cite{GH17}).\par

The above decomposition of the kernel as a vGBS group allows for a very explicit description of all the fibrations of a RAAG $\raag \Gamma$ when the underlying graph $\Gamma$ is a block graph; in this case edge groups are trivial, and  we can obtain a  decomposition of the Artin kernel as a free product of free abelian groups, whose rank can explicitly be computed, leading to the following rank formula (see Theorem \ref{thm:rankformula}, and compare the rank formula in \cite[Theorem 1.1,6.1]{CL16}). Here $\block \Gamma$ is the set of blocks (i.e. biconnected components, or equivalently maximal cliques) of the block graph $\Gamma$, and $\bldeg_\Gamma(v)$ is the number of blocks containing a vertex $v$.

\begin{theorem}
Let $\Gamma$ be a connected block graph, and $f:A_\Gamma\to \mathbb Z$ be non-zero on cut vertices. Then

$$\rk( \ker(f))= 1+  \sum_{B\in \block \Gamma} [f(\raag \Gamma):f(\raag{B})](|\vv B|-2)  +\sum_{v\in \vv \Gamma}\dfrac{\left( \bldeg_\Gamma(v)-1 \right)|f(v)|}{[\zz:f(\raag{\Gamma})]} $$
\end{theorem}

As an application of the explicit control that this formula provides, we show that when $\Gamma$ is a block graph, the RAAG $\raag \Gamma$ admits fibrations with kernels of arbitrarily large rank (see Corollary \ref{cor:unbounded_fibers}), and we characterize those with kernel of minimal rank (see Corollary \ref{cor:minimal_fibers}). Moreover in Corollary \ref{cor:bounded_div} we show how to construct sequences of tame fibrations
$$1\to K_n \to \raag \Gamma \overset{f_n}{\to} \zz \to 1$$
whose kernels have constant finite rank, but such that the sequence $f_n$ converges to a homomorphism $f_\infty :\raag \Gamma \to \zz$ whose kernel is wild (i.e. surjects onto $\ff_\infty$). Here the convergence is best understood to be up to rescaling by (positive) real numbers, operation that of course does not modify the kernel; in this sense this shows an explicit way in which the function that records the rank of the kernel of a fibration is not a proper function on the BNS invariant of $\raag \Gamma$.

\vspace{.5cm}

\textbf{Acknowledgements}: Part of this project was supported by the Research Foundation Flanders (Project G.0F93.17N), and was developed at KU Leuven (Kulak); part of this project was developed in the framework of UROP at FSU.
We would like to thank the organizers of the 2017-2018 Warwick EPSRC Symposium on Geometry, Topology and Dynamics in Low Dimensions for organizing very fruitful workshops, as well as Conchita Mart\'inez P\'erez and Armando Martino for useful conversations. We also thank Ian Leary for suggesting a change in the terminology that was used in a preliminary version of this paper.

\section{Background on right-angled Artin groups}
Let $\Gamma$ be a graph (by which we mean a finite 1-dimensional simplicial complex). We denote by $\vv{\Gamma}$ the set of its vertices and by $\ee{\Gamma}$ the set of its edges. The right-angled Artin group (RAAG in the following) associated to $\Gamma$ is the group $\raag{\Gamma}$ generated by the vertices of $\Gamma$ and such that two generators commute if and only if they are joined by an edge, i.e. $\raag{\Gamma}$ has the following presentation

$$\raag{\Gamma} = \langle x \in \vv{\Gamma} \ | \ [x,y] \in \ee{\Gamma} \rangle$$

The RAAG associated to the empty graph is defined to be the trivial group. The RAAG associated to a clique (i.e. a complete graph) on $n$ vertices is the free abelian group $\zz^n$. The RAAG associated to a totally disconnected graph on $n$ vertices is the free group $\ff_n$.

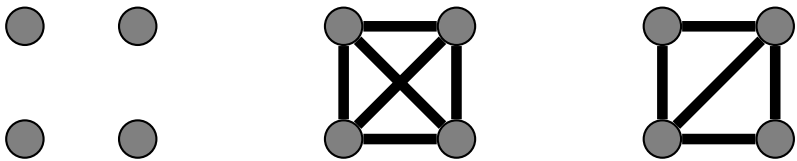
\begin{figure}[h]
  \begin{minipage}[c]{0.33\linewidth}
    \centering
    
    \begin{tikzpicture}[thick]
    \tikzstyle{every node}=[circle, draw, fill=black!50,
                        inner sep=5pt, minimum width=2pt]
    
    
        \node (n0)  at (0,0) {};
        \node (n1)  at (0,1.5) {};
        \node (n2)  at (1.5,0) {};
        \node (n3)  at (1.5,1.5) {};
        
        
    \end{tikzpicture}
  \end{minipage}
  \begin{minipage}[c]{0.34\linewidth}
    \centering
    
    \begin{tikzpicture}[thick]
    \tikzstyle{every node}=[circle, draw, fill=black!50,
                        inner sep=5pt, minimum width=2pt]
    
    
        \node (n0)  at (0,0) {};
        \node (n1)  at (0,1.5) {};
        \node (n2)  at (1.5,0) {};
        \node (n3)  at (1.5,1.5) {};
    
    \path[every node/.style={font=\sffamily},line width=4pt]
    
        (n0) edge node {} (n1)
        (n0) edge node {} (n2)
        (n0) edge node {} (n3)
        (n1) edge node {} (n2)
        (n1) edge node {} (n3)
        (n2) edge node {} (n3);
        
    \end{tikzpicture}
  \end{minipage}
  \begin{minipage}[c]{.33\linewidth}
    \centering
    
    \begin{tikzpicture}[thick]
    \tikzstyle{every node}=[circle, draw, fill=black!50,
                        inner sep=5pt, minimum width=2pt]
    
    
        \node (n0)  at (0,0) {};
        \node (n1)  at (0,1.5) {};
        \node (n2)  at (1.5,0) {};
        \node (n3)  at (1.5,1.5) {};
    
    \path[every node/.style={font=\sffamily},line width=4pt]
    
        (n0) edge node {} (n1)
        (n0) edge node {} (n2)
        (n0) edge node {} (n3)
        (n1) edge node {} (n3)
        (n2) edge node {} (n3);
        
    \end{tikzpicture}
  \end{minipage}
    \caption{Graphs corresponding to $\ff_4$, $\zz^4$, and an amalgamated product of two copies of $\zz^3$ over $\zz^2$.}
\end{figure}

More general RAAGs interpolate between these two examples, and can be quite rich despite their simple definition; for instance when $\Gamma$ is a tree $\raag \Gamma$ is the fundamental group of a compact 3-manifold (with boundary) fibering over $S^1$ (see \cite{D87}).\par 
It is well-known (see \cite{D87b}) that the graph completely determines the group, in the sense that $\raag \Gamma$ is isomorphic to $\raag \Lambda$ if and only if the graphs $\Gamma$ and $\Lambda$ are isomorphic. Moreover given any full subgraph $\Lambda\subseteq \Gamma$ we have that the RAAG $\raag \Lambda$ embeds in a natural way as a subgroup of $\raag \Gamma$.

\begin{remark}\label{rmk:maptoabel}
Notice that any map of $\raag \Gamma$ into a group is completely determined by what it does on vertices. Moreover, since all the relations are just commutators, for any abelian group $A$ any map $f_0:\vv \Gamma \to A$ uniquely extends to a well-defined group homomorphism $f:\raag \Gamma\to A$.
\end{remark}

Bestvina and Brady in \cite{BB97} considered the map $f:\raag \Gamma\to \zz$ obtained by sending every vertex to $1\in \zz$. The kernel of this map is known as the Bestvina-Brady subgroup of $\raag \Gamma$. In this paper we consider the following generalization. 

\begin{definition}
An \textit{Artin kernel} is a normal subgroup $K$ of $\raag \Gamma$ such that the quotient $\raag \Gamma / K$ is infinite cyclic.
\end{definition}

A practical way to construct such subgroups is of course to consider kernels of more general maps that send the vertices of $\Gamma$ to more general collections of numbers (see Section \ref{sec:characters} below).  Finiteness properties of Artin kernels have been studied by Bux and Gonzalez in \cite{BG99}, where they have extended the original results by Bestvina and Brady; further algebraic properties have been studied by Papadima and Suciu in \cite{PS09}. \par

Notice that, unlike right-angled Artin groups, Artin kernels are very far away from being classified by their underlying graph: for example the Bestvina-Brady group on a tree is always a free group generated by the edges. Moreover, while certain Artin kernels are themselves RAAGs, this does not hold for general Artin kernels, for instance because they can fail to be finitely presented (see \cite{BG99}, and for an example even in the finitely presented case see \cite[Example 2.8]{PS07}).\par
We would also like to warn the reader that while Artin kernels provide a natural generalization of the classical Bestvina-Brady groups, the name ``generalized Bestvina-Brady groups'' is nowadays mostly being used to refer to the groups introduced by Leary in \cite{L18}; these groups are not subgroups of RAAGs in general.

\subsection{Graphical splittings for RAAGs}

We will be interested in understanding the structure of Artin kernels induced by decompositions of the underlying graph $\Gamma$.
We start by reviewing known results about splittings of RAAGs induced by decompositions of the underlying graph.\par

Recall that a group $G$ splits over a subgroup $H\neq G$ if it can be realized as the free product of two subgroups $A,B\neq G$  amalgamated over $H$ (i.e. $G\cong A \ast_H B$), or as an HNN extension of a subgroup $A\neq G$ over $H$ (i.e. $G\cong A \ast_H$).
For easy ways to obtain splittings of a RAAG one can notice the following: if the graph is disconnected then $\raag \Gamma$ splits as a free product of the RAAGs over the connected components, and if it is complete then $\raag \Gamma$ is free abelian, hence splits as an HNN extension over a codimension-1 subgroup.
Less trivial splittings can be induced by decompositions of the graph, as follows.

\begin{definition}\label{def:separating}
For a subgraph $\Lambda \subseteq \Gamma$ we denote by $\Gamma \setminus \Lambda$ the full subgraph of $\Gamma$ generated by $\vv \Gamma \setminus \vv \Lambda$, and call it the \textit{complement} of $\Lambda$ in $\Gamma$. We say $\Lambda$ is \textit{separating} when its complement is disconnected. A \textit{cut vertex} is a separating vertex.
\end{definition}

\begin{definition}\label{def:splitting}
Let $\Gamma_1, \Gamma_2, \Gamma_3$ be pairwise distinct non-empty  full subgraphs of $\Gamma$ such that $\Gamma=\Gamma_1\cup \Gamma_2$ and $\Gamma_1 \cap \Gamma_2 = \Gamma_3$.
We say that $(\Gamma_1,\Gamma_2,\Gamma_3)$ is a \textit{splitting} of $\Gamma$. If $\Gamma_3$ is connected, then we say $(\Gamma_1,\Gamma_2,\Gamma_3)$ is a \textit{connected splitting} of $\Gamma$.
\end{definition}

For the sake of precision let us emphasize that here the condition $\Gamma=\Gamma_1\cup \Gamma_2$ means that both $\vv \Gamma=\vv {\Gamma_1}\cup \vv {\Gamma_2} $ and  $\ee \Gamma=\ee {\Gamma_1}\cup \ee {\Gamma_2}$, and the condition $\Gamma_3=\Gamma_1\cap \Gamma_2$ means that both $\vv {\Gamma_3}=\vv {\Gamma_1}\cap \vv {\Gamma_2} $ and  $\ee { \Gamma_3}=\ee {\Gamma_1}\cap \ee {\Gamma_2}$. In particular we have the following.

\begin{lemma}\label{lem:graph_props}
 Let $( \Gamma_1, \Gamma_2, \Gamma_3)$ be a splitting of $\Gamma$. Then the following hold.
\begin{enumerate}
    \item \label{item:gp1} There is no edge between $\vv {\Gamma_1}\setminus \vv{ \Gamma_3}$ and $\vv{\Gamma_2}\setminus \vv{ \Gamma_3}$.
    
    \item \label{item:gp2} $\Gamma_3$ is separating.
    
    \item \label{item:gp3}  If $\Gamma$ and $\Gamma_3$ are connected, then  $\Gamma_1$ and $\Gamma_2$ are connected too.
\end{enumerate}
\end{lemma}
\proof
Statement (\ref{item:gp1}) follows from the fact that $\Gamma=\Gamma_1\cup \Gamma_2$ and  $\Gamma_1 \cap \Gamma_2 = \Gamma_3$, and  statement (\ref{item:gp2}) follows from (\ref{item:gp1}).
To prove (\ref{item:gp3}), notice that since $\Gamma$ is connected, we can join any vertex of $\Gamma_1$ to a vertex in $\Gamma_3$ by an edge-path; moreover by (\ref{item:gp1}) this edge-path can be chosen to lie completely in $\Gamma_1$. Then connectedness of  $\Gamma_3$ implies that any two vertices in $\Gamma_1$ can be joined by an edge-path entirely contained in $\Gamma_1$; the same applies to $\Gamma_2$.
\endproof

It is easy to check directly from the RAAG presentation that a splitting of $\Gamma$ induces a splitting of $\raag \Gamma$, in the following sense.
\begin{lemma}\label{lem:split}
If $(\Gamma_1,\Gamma_2,\Gamma_3)$ split $\Gamma$, then $\raag \Gamma=\raag{\Gamma_1}\underset{\raag{\Gamma_3}}{\ast}\raag{\Gamma_2}$.
\end{lemma}

For abelian splittings (i.e. splittings over abelian subgroups) actually more is true, as shown by Clay for splitting over $\zz$ (in \cite{C14}), and Groves and Hull in the general case (see \cite{GH17}). Notice that an abelian subgroup of a RAAG is necessarily torsion-free.

\begin{proposition}\label{prop:easy_split} 
$\raag \Gamma$ splits over an abelian group if and only if one of the following occurs:
\begin{enumerate}
  \item $\Gamma$ is disconnected, 
  \item $\Gamma$ is complete, 
  \item $\Gamma$ has a separating clique.
 \end{enumerate}
\end{proposition}

In particular one can see that when $\Gamma$ is disconnected, $\raag \Gamma$ splits as a free product (indeed Lemma \ref{lem:split} holds even when $\Gamma_3=\varnothing$), and the same holds for any Artin kernel defined on the same graph.
On the other hand when $\Gamma$ is complete both $\raag \Gamma$ and all its Artin kernels are free abelian groups; as such they all split as HNN extensions over codimension-1 subgroups.
Therefore in the rest of the paper we will be interested in the case in which $\Gamma$ is connected and not complete, hence contains a  separating subgraph.

\subsection{Characters and BNS invariant for RAAGs}\label{sec:characters}
Normal subgroups with infinite cyclic quotients can be conveniently defined as kernels of real characters, by embedding the quotient in $\rr$.
In this section we develop this point of view by reviewing a convenient space of real characters introduced by Bieri-Neumann-Strebel in \cite{BNS87} for general finitely presented groups, and then studied by Meier and VanWyk in the specific case of RAAGs in \cite{MV95}. We are going to use this to organize the collection of Artin kernels of a fixed RAAG, and to have a meaningful notion of convergence.\par

Let $\raag \Gamma$ be a RAAG defined by a connected graph $\Gamma$. The space of all real characters (i.e. group homomorphisms into $\rr$) can be easily parametrized by  recording the value of a character on the vertices of $\Gamma$, i.e.
$$\operatorname{Hom}(\raag \Gamma,\rr)\overset{\cong}{\to} \rr^{|\vv \Gamma|}, \quad f\mapsto (f(v_1),\dots,f(v_n))$$
where $\vv \Gamma = \{v_1,\dots,v_n\}$;  notice this is as isomorphism by Remark \ref{rmk:maptoabel}.

We explicitly notice that if $f:\raag \Gamma\to \rr$ is a character whose image is a discrete non-trivial (equivalently, infinite cyclic) subgroup of $\rr$, then $\ker(f)$ is an  Artin kernel; conversely if $K\subset \raag \Gamma$ is an Artin kernel, then $\raag \Gamma / K \cong \zz$, and any isomorphism provides a character of $\raag \Gamma$.
In other words, characters with discrete image correspond to Artin kernels.
We find it useful to introduced the following terminology.

\begin{definition}
We say a character $f:\raag \Gamma \to \rr$ is respectively \textit{integral}, \textit{rational} or \textit{discrete} if its image is an infinite cyclic subgroup of $\zz, \qq$ or $\rr$.
\end{definition}

Rescaling a character does not change its kernel, and a character with discrete image can always be rescaled to be integral, so that Artin kernels can actually be thought as primitive integral vectors in $\rr^{|\vv \Gamma|}$ through the above isomorphism.
Notice that the zero character is not discrete according to this definition. To avoid useless complications, unless otherwise specified, in the following we will implicitly assume that $f$ is not the zero character. \par
As we are interested in characters only up to rescaling, we also find it convenient to introduce the \textit{character sphere }
$$\chars {\raag \Gamma}=\left(\operatorname{Hom}(\raag \Gamma,\rr)\setminus \{0\} \right) / \rr^+$$
and to conflate a character and its class in $\chars {\raag \Gamma}$, whenever no confusion arises. We  will also extend the above terminology accordingly. 
We remark that here we are taking the quotient only with respect to positive scalars (as opposed to non-zero scalar) just for the sake of consistency with the literature. \par

Meier and VanWyk in \cite{MV95} have obtained a simple condition for a discrete character of a RAAG to have a finitely generated kernel (notice they call rational what we call discrete here).

\begin{definition}
Given a character $f:\raag \Gamma \to \rr$ we define
\begin{itemize}
    \item the \textit{living subgraph} of $f$ as the full subgraph $\mathcal L (f)$ of $\Gamma$ generated by $\{v\in \vv \Gamma \ | \ f(v)\neq 0  \}$;
    
    \item the \textit{dead subgraph} of $f$ as the full subgraph $\mathcal D(f)$ of $\Gamma$ generated by $\{v\in \vv \Gamma \ | \ f(v)=0  \}$.
\end{itemize} 
\end{definition}
Notice these graphs depend only on the equivalence class of $f$. Moreover when $\Gamma$ is connected we have that $\mathcal L(f)$ is connected if and only if $\mathcal D(f)$ is not separating.
\begin{definition}
A subgraph $\Lambda$  of $\Gamma$ is \textit{dominating} if every vertex of $\Gamma$ is adjacent to some vertex of $\Lambda$, i.e. $\Gamma$ coincides with the star of $\Lambda$.
\end{definition}

As far as our discussion is concerned,  the  main result by Meier and VanWyk  is the following.
\begin{theorem}[\textrm{\cite[Theorem 6.1]{MV95}}] \label{thm:living_subgraph}
Let $f:\raag \Gamma\to \rr$ be discrete. Then $\ker(f)$ is finitely generated if and only if $ \mathcal L (f)$ is connected and dominating. 
\end{theorem}

The \textit{Bieri-Neumann-Strebel invariant} (in the following \textit{BNS invariant}) of a finitely generated group was introduced in \cite{BNS87} as a device to capture information about finiteness properties of normal subgroups with abelian quotients, and therefore it is  a convenient way to organize characters with finitely generated kernels. It is an open subset of the character sphere 
$$\Sigma^1({\raag \Gamma})\subseteq  \chars {\raag \Gamma}$$
which can be thought as a generalization of Thurston's polyhedron for 3-manifold groups (see \cite{T86}). The original definition is not the most useful in our setting, and we prefer to use the following characterization (see \cite[\S 4]{BNS87} or \cite[Corollary 1.2]{MV95}) as a working definition of what it means for a discrete character to belong to the BNS invariant.

\begin{theorem}[Bieri-Neumann-Strebel]\label{thm:BNS}
Let $f:\raag \Gamma \to \rr$ be discrete. Then $\ker(f)$ is finitely generated if and only if both $f$ and $-f \in \Sigma^1(\raag \Gamma)$.
\end{theorem}

To sum up we have the following characterization of discrete characters with finitely generated kernels, i.e. finitely generated Artin kernels of a fixed RAAG.

\begin{corollary}\label{cor:MV+BNS}
Let $f:\raag \Gamma \to \rr$ be discrete. Then the following are equivalent:
\begin{enumerate}
\item \label{item:fg} $\ker(f)$ is finitely generated,
\item \label{item:living} $\mathcal L(f)$ is connected and dominating,
\item \label{item:both} $f,-f \in \Sigma^1(\raag \Gamma)$ ,
\item \label{item:S1} $f \in \Sigma^1(\raag \Gamma)$.
\end{enumerate}
\end{corollary}
\proof 
The equivalence of (\ref{item:fg}) and (\ref{item:living}) is Theorem \ref{thm:living_subgraph}, and the equivalence of (\ref{item:fg}) and (\ref{item:both}) is Theorem \ref{thm:BNS}. The equivalence of (\ref{item:both}) and   (\ref{item:S1}) follows from the fact that  $\Sigma^1(\raag \Gamma)$ is invariant under automorphisms of $\raag \Gamma$, and the fact that any RAAG has an involutive automorphism which sends every generator to its inverse; the action of such an automorphism on characters is of course $f\mapsto -f$.
\endproof

\section{Splittings of Artin kernels}

We now turn to the study of the finiteness properties of the Artin kernels of a fixed RAAGs in terms of the behavior of a defining  character on  separating subgraphs.
Recall that in order to avoid trivial  redundancies, by default we assume that $f:\raag \Gamma \to \rr$ is not the zero character, and that  $\Gamma$ is a finite connected non-complete graph, so that it admits non-trivial splittings. \par
After proving a general dichotomy relative to a single splitting in \S \ref{sec:dichotomy}, in \S \ref{sec:chordal} we focus on   graphs that can be iteratively decomposed into pieces which are easily understood, such as chordal graphs. In the last section \S \ref{sec:block_graphs} we obtain an explicit formula for the rank of the Artin kernels and discuss some features of the BNS invariant of $\raag \Gamma$ when $\Gamma$ is a block graph.

\subsection{Splitting dichotomy}\label{sec:dichotomy}

Here we prove a general splitting result along a separating subgraph, which will later be used as a key tool for more specific applications. More precisely we propose an extension to Artin kernels of the main Proposition in \cite{D87}, which is about the classical Bestvina-Brady groups. This is inspired by an analogous dichotomy obtained by Cashen and Levitt in \cite[Theorem 4.4]{CL16} in the context of free-by-cyclic groups. It should be noted that Bestvina-Brady groups are free-by-cyclic precisely when the underlying graph is a tree. The main argument is based on the following  orbit count lemma (which appears in \cite[Lemma 4.5]{CL16}), which we will apply to the action of Artin kernels on Bass-Serre trees associated to splittings of the underlying graph. 

\begin{lemma}\label{lem:orbitcount}
Let a group $G$ act transitively on a set $X$; let $K\subseteq G$ be a normal subgroup and $f:G\to G/K$ be the quotient map.
For $x\in X$ let $G_x$ denote the stabilizer of $x$ for the action $G\curvearrowright X$.
Then the number of $K$-orbits in $X$ is equal to the index $[G/K:f(G_x)]$.
\end{lemma}

Here and in the following, for any subgraph $\Lambda$ of $\Gamma$ we let $\ind_{\Gamma,\Lambda}=[f(\raag \Gamma):f(\raag{\Lambda})]$ be the index of $f(\raag{\Lambda})$ in $f(\raag{\Gamma})$. We begin with a statement that does not require the splitting to be connected. We also remark that by a finite graph of groups we just mean a graph of groups whose underlying graph is finite, but in which groups are allowed to be not finitely generated.

\begin{theorem}\label{thm:generaldichotomy}
Let $(\Gamma_1,\Gamma_2,\Gamma_3)$ be a splitting of $\Gamma$, let $f:\raag \Gamma \to \rr$ be discrete, and let $f_k:\raag {\Gamma_k} \to \rr$ be its restriction for $k=1,2,3$.
Then $\ker (f)$ splits as a graph of groups $\mathcal G$ with 
    \begin{itemize}
    \item $\ind_{\Gamma,\Gamma_k}$ vertices with group isomorphic to $\ker(f_k)$ for $k=1,2$,
    \item $\ind_{\Gamma,\Gamma_3}$ edges with group isomorphic to $\ker(f_3)$.
    \end{itemize}
In particular 
\begin{enumerate}
\item \label{item:finitegog} if  $f_3\neq 0$, then $\ker (f)$ splits as a finite graph of groups $\mathcal G$,
\item \label{item:surjectFinf} if  $f_3= 0$, then $\ker (f)$ surjects onto $\ff_\infty$. 
\end{enumerate}
\end{theorem}
\proof
By Lemma \ref{lem:split}, we know that the splitting of $\Gamma$ induces a splitting of $\raag \Gamma$ as

$$\raag \Gamma=\raag{\Gamma_1}\underset{\raag{\Gamma_3}}{*}\raag{\Gamma_2}$$
Let $\mathcal T$ be the Bass-Serre tree associated to this splitting.
The group $\raag \Gamma$ acts on $\mathcal T$, and we restrict the action to the normal subgroup $\ker(f)$.
Let $x\in \mathcal{T}$ be a vertex corresponding to a coset of $\raag{\Gamma_1}$.
We have that $\raag \Gamma$ acts transitively on the orbit $\raag \Gamma . x$, with vertex stabilizers conjugate to 
$\raag{\Gamma_1}$; by Lemma \ref{lem:orbitcount} this orbit splits in $\ind_{\Gamma,\Gamma_1}$ orbits for the action of 
$\ker(f)$, and vertex stabilizers are conjugate to $\ker(f)\cap \raag{\Gamma_1}=\ker (f_1)$. A similar argument applies to the other type of vertices, and to the edges. Then the quotient $\mathcal G=\mathcal T/\ker(f)$, labelled with the appropriate stabilizers, provides the desired graph of groups decomposition. Notice that this graph of groups is bipartite and has no loops. We are left to proving the dichotomy.
\begin{enumerate}
    \item If $f_3\neq 0$, then of course also $f_k\neq 0$ for $k=1,2$, hence $\ind_{\Gamma,\Gamma_k}$ is finite for $k=1,2,3$, hence from the first part of the statement we get that $\mathcal G$ is finite.
    \item If $f_3= 0$, then $\ind_{\Gamma,\Gamma_3}=\infty$, and the graph of groups has infinitely many edges. We distinguish two cases.
    \begin{enumerate}
        \item If $f_k\neq 0$ for $k=1,2$, then $\mathcal G$ has finitely many vertices, hence infinitely many bigons; in particular the underlying graph has  fundamental group $\ff_\infty$. Killing all the vertex groups provides a surjection to this infinitely generated group.
        \item If $f_1= 0$, then necessarily $f(\raag{ \Gamma_2})=f(\raag \Gamma)\neq 0$ because $f$ is not the trivial character. 
        In particular $\mathcal G$ has infinitely many vertices corresponding to $\Gamma_1$, and a single vertex corresponding to $\Gamma_2$, i.e. it is a star with infinitely many edges. Notice that in this sub-case we have $\ker(f_1)=\raag{\Gamma_1}$, as well as $\ker(f_3)=\raag{\Gamma_3}$.
        Since $\Gamma_3\neq \Gamma_1$, each of the terminal vertex groups admits a non-trivial map to $\zz$ which kills (at least) the incoming edge group. Combining these maps one obtains the desired surjection to an infinitely generated free group.
        The case $f_2=0$ is totally analogous to the case $f_1=0$.
    \end{enumerate}
\end{enumerate}

\endproof

\begin{figure}[ht]
    \centering
  \begin{minipage}[c]{0.5\linewidth}
    \centering
    
    \begin{tikzpicture}[thick]
    \tikzstyle{every node}=[circle, draw, fill=black!50,
                        inner sep=5pt, minimum width=2pt]
    
    
        \node[fill=blue, label=left:{$1$}] (n0)  at (-1,0) {};
        \node[fill=red, label=above:{$1$}] (n1) at (0,1) {};
        \node[fill=green!70!black, label=right:{$0$}] (n2) at (1,0) {};
        \node[fill=red, label=below:{$1$}] (n3) at (0,-1) {};
    
    \path[every node/.style={font=\sffamily},line width=4pt]
    
        (n0) edge[blue] node {} (n1)
        (n1) edge[green!70!black] node {} (n2)
        (n2) edge[green!70!black] node {} (n3)
        (n3) edge[blue] node {} (n0);
        
    \end{tikzpicture}
  \end{minipage}
  \begin{minipage}[c]{0.5\linewidth}
    \centering
    
    \begin{tikzpicture}[thick]
    \tikzstyle{every node}=[circle, draw, fill=black!50,
                        inner sep=5pt, minimum width=2pt]
    
    
        \node[fill=blue, label=left:{$1$}] (n0)  at (-1,0) {};
        \node[fill=red, label=above:{$1$}] (n1) at (0,1) {};
        \node[fill=green!70!black, label=right:{$1$}] (n2) at (1,0) {};
        \node[fill=red, label=below:{$0$}] (n3) at (0,-1) {};
    
    \path[every node/.style={font=\sffamily},line width=4pt]
    
        (n0) edge[blue] node {} (n1)
        (n1) edge[green!70!black] node {} (n2)
        (n2) edge[green!70!black] node {} (n3)
        (n3) edge[blue] node {} (n0);
        
    \end{tikzpicture}
  \end{minipage}
  \begin{minipage}[c]{0.5\linewidth}
    \centering
    
    \begin{tikzpicture}[thick]
    \tikzstyle{every node}=[circle, draw, fill=black!50,
                        inner sep=5pt, minimum width=2pt]
    
    
        \node[fill=blue, label=left:{$1$}] (n0)  at (-1,0) {};
        \node[fill=red, label=above:{$0$}] (n1) at (0,1) {};
        \node[fill=green!70!black, label=right:{$1$}] (n2) at (1,0) {};
        \node[fill=red, label=below:{$0$}] (n3) at (0,-1) {};
    
    \path[every node/.style={font=\sffamily},line width=4pt]
    
        (n0) edge[blue] node {} (n1)
        (n1) edge[green!70!black] node {} (n2)
        (n2) edge[green!70!black] node {} (n3)
        (n3) edge[blue] node {} (n0);
        
    \end{tikzpicture}
  \end{minipage}
  \begin{minipage}[c]{0.5\linewidth}
    \centering
    
    \begin{tikzpicture}[thick]
    \tikzstyle{every node}=[circle, draw, fill=black!50,
                        inner sep=5pt, minimum width=2pt]
    
    
        \node[fill=red, label=left:{$1$}] (n0)  at (-1,0) {};
        \node[fill=blue, label=above:{$0$}] (n1) at (0,1) {};
        \node[fill=red, label=right:{$1$}] (n2) at (1,0) {};
        \node[fill=green!70!black, label=below:{$0$}] (n3) at (0,-1) {};
    
    \path[every node/.style={font=\sffamily},line width=4pt]
    
        (n0) edge[blue] node {} (n1)
        (n1) edge[blue] node {} (n2)
        (n2) edge[green!70!black] node {} (n3)
        (n3) edge[green!70!black] node {} (n0);
        
    \end{tikzpicture}
  \end{minipage}
    
    \caption{In the first three graphs the splitting subgraph is given by the union of the top and the bottom vertices, while in the last one it is given by the union of the left and the right vertices.}
    \label{fig:ex_different_splittings}
\end{figure}
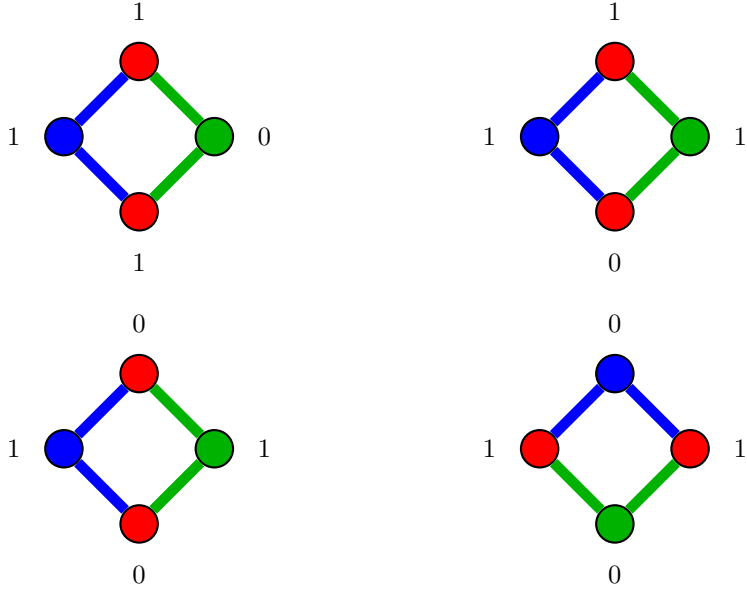
\begin{remark}\label{rmk:general_dic_caveat}
Notice that the separating graph does not have to be connected, and indeed in case (\ref{item:finitegog}) it is possible for either edge groups or vertex groups (or both) to be infinitely generated. Moreover it is possible for a fixed character $f$ to behave as in case (\ref{item:finitegog}) with respect to one splitting and as in case  (\ref{item:surjectFinf}) with respect to another one, i.e. $\ker(f)$ can both surject to $\ff_\infty$ and split as a finite graph of groups (some vertex groups will then surject to $\ff_\infty$).
Concrete examples of these phenomena can be constructed on any cycle of at least four vertices, see Figure \ref{fig:ex_different_splittings} for some concrete examples. We will see that this does not happen for chordal graphs in \S \ref{sec:chordal}.

\end{remark}

It is known by \cite[Theorem 6.1]{MV95} that if $f=0$  on a separating subgraph then $\ker(f)$ is not finitely generated (compare Theorem \ref{thm:living_subgraph}).
The above proof provides an independent proof of this fact, and actually shows more: as soon as $\ker(f)$ is not finitely generated it even surjects onto $\ff_\infty$. 
The next statement is a sort of converse to case (\ref{item:surjectFinf}) of Theorem \ref{thm:generaldichotomy}. Recall that the dead subgraph $\mathcal D(f)$ of $f$ is defined to be the full subgraph generated by the vertices that $f$ sends to $0$.

\begin{corollary}[Wild Artin kernels]\label{cor:wild}
Let $f:\raag \Gamma \to \rr$ be discrete. Then the following are equivalent:
\begin{enumerate}
    \item \label{item:fkillssep} $\exists \Lambda \subseteq \Gamma$ separating subgraph such that $ \Lambda \subseteq \mathcal D(f)$ (i.e. $f=0$ on $\Lambda$).
    
    \item \label{item:kerfsurjects} $\ker(f)$ surjects onto $\ff_\infty$.
    
    \item \label{item:kerfnotfg} $\ker(f)$ is not finitely generated.
\end{enumerate}

\end{corollary}
\proof
The fact that (\ref{item:fkillssep}) implies (\ref{item:kerfsurjects}) follows from Theorem \ref{thm:generaldichotomy}, case (\ref{item:surjectFinf}); of course (\ref{item:kerfsurjects}) implies (\ref{item:kerfnotfg}). So we only need to show that  (\ref{item:kerfnotfg}) implies (\ref{item:fkillssep}).
Suppose $\ker(f)$  is not finitely generated. Then by Theorem \ref{thm:living_subgraph} we know that $\mathcal{L}(f)$ is either not connected or not dominating. We distinguish two cases.
\begin{enumerate}
    \item If $\mathcal{L}(f)$ is not connected, then  $\mathcal D(f)$ is separating.
    
    \item If $\mathcal{L}(f)$ is not dominating, then there is a vertex $v\in \vv \Gamma$ whose link entirely consists of vertices on which $f$ is zero, i.e. $lk(v)\subseteq \mathcal D(f)$. By definition $lk(v)$ is separating. 
\end{enumerate}
\endproof

In particular as soon as $\Gamma$ is not complete, it has a separating subgraph, hence there is a wild Artin kernel, i.e. a normal subgroup $ K \subseteq \raag \Gamma$ such that $\raag \Gamma / K$ is infinite cyclic and $K$ surjects to $\ff_\infty$.\par

So far the separating graphs have been allowed to be disconnected. On the other hand in order to obtain stronger  finiteness properties for the tame case, we need to work with connected splittings (see Definition~\ref{def:splitting}). Recall that by Lemma \ref{lem:graph_props} all the subgraphs in a connected splitting of a connected graph are connected.
The next statement should be compared to the general treatment of BNS invariants for graphs of groups in \cite[\S 2]{CL16}; in this regard, it should also be noted that non-abelian RAAGs are not slender because they contain $\ff_2$ (recall that a group is said to be slender when every subgroup is finitely generated).

\begin{lemma}\label{lem:BNSoversplit} 
Let $(\Gamma_1,\Gamma_2,\Gamma_3)$ be a splitting of $\Gamma$, let $f:\raag \Gamma \to \rr$ be discrete, and let $f_k:\raag {\Gamma_k} \to \rr$ be its restriction for $k=1,2,3$. Then the following hold.
\begin{enumerate}
    \item \label{item:restriction_discrete} $f_k$ is zero or discrete.
    \item \label{item:restriction_f3nonzero} If $\ker(f)$ is finitely generated, then $f_3$ is non-zero (and discrete).
    
    \item \label{item:restriction_fg} If $\ker(f)$ is finitely generated and $\ker(f_3)$ is finitely generated, then $\ker(f_k)$ is finitely generated for $k=1,2$.
    
\end{enumerate}
\end{lemma}
\proof

The restriction of a discrete character is clearly a discrete character as soon as it is not the zero character, so (\ref{item:restriction_discrete}) is trivial. To prove (\ref{item:restriction_f3nonzero}), assume by contradiction that $f_3$ is identically $0$ on $\Gamma_3$. Notice $\Gamma_3$ is a separating subgraph, contained in $\mathcal D(f)$. By Corollary~\ref{cor:wild} we get that $\ker(f)$ can not be finitely generated, which is a contradiction.\par

Finally let us prove (\ref{item:restriction_fg}). Notice that the hypotheses imply that $\Gamma$ and $\Gamma_3$ are connected. By Lemma \ref{lem:graph_props} we get that $\Gamma_1,\Gamma_2$ are connected too.
Moreover by Corollary \ref{cor:MV+BNS} we get that $\mathcal L(f)$ (respectively $\mathcal L(f_3)$) is connected and dominating as a subgraph of $\Gamma$ (respectively $\Gamma_3$). Notice that $\mathcal L(f_k)=\mathcal L(f) \cap \Gamma_k$.
So, focusing on $\Gamma_1$ (the same proof works for $\Gamma_2$), we need to show that $\mathcal L(f_1)$ is connected and dominating as a subgraph of $\Gamma_1$. \par

Since $\mathcal L (f)$ is connected, any two vertices $x,y$ in $\mathcal L(f_1)$ can be joined by an edge-path contained in $\mathcal L (f)$. We claim that such a path can be chosen to be also contained in $\Gamma_1$. By contradiction assume that any edge-path from $x$ to $y$ in 
$\mathcal L (f)$  makes an excursion into $\Gamma \setminus \Gamma_1=\Gamma_2\setminus \Gamma_3$; since $\Gamma_3$ is separating every such excursion starts and ends with a vertex on $\mathcal L(f_3)$. But  $\mathcal L (f_3)$ is connected, so there is a shortcut entirely contained in $\mathcal L(f_3)\subseteq \mathcal L (f_1)$, leading to a contradiction. \par

To see that $\mathcal L(f_1)$ is  dominating as a subgraph of $\Gamma_1$ we can argue as follows: if $x\in \vv{\Gamma_3}$ then we use that $\mathcal L(f_3)$ is dominating in $\Gamma_3$ to find a vertex in $\mathcal L(f_3)$ at distance 1 from $x$. If $x\in \vv{\Gamma_1}\setminus \vv{\Gamma_3}$, then since $\mathcal L(f)$ is dominating in $\Gamma$ we find a vertex of $\mathcal L(f)$ at distance 1 from $x$; notice that by Lemma \ref{lem:graph_props} this vertex can not live in $\vv{\Gamma_2}\setminus \vv{ \Gamma_3}$.
\endproof

We then deduce the following statement. Notice its hypothesis are satisfied trivially when $f$ does not vanish on any vertex, i.e. $\mathcal{D}(f)= \varnothing$.
\begin{corollary}[Tame Artin kernels]\label{cor:tame}
Let $(\Gamma_1,\Gamma_2,\Gamma_3)$ be a connected splitting of $\Gamma$, and let $f:\raag \Gamma \to \rr$ be a discrete character such that $\ker(f)$ and $\ker(f_3)$ are finitely generated. 
Then $\ker (f)$ splits as a finite graph of groups $\mathcal G$ as in Theorem \ref{thm:generaldichotomy} with finitely generated edge and vertex groups.
\end{corollary}

\proof
From Theorem \ref{thm:generaldichotomy} we know that $\ker(f)$ splits as a finite graph of groups $\mathcal G$ in which vertex groups are isormorphic to $\ker(f_k)$ for $k=1,2$ and edge groups are isomorphic to $\ker(f_3)$. So the statement boils down to showing that $\ker(f_k)$ is finitely generated for $k=1,2,3$, which follows from  Lemma \ref{lem:BNSoversplit}.
\endproof

\subsection{Chordal graphs}\label{sec:chordal}

In the previous section we have identified conditions under which a discrete character displays a wild behavior (i.e. surjects to $\ff_\infty$) and conditions under which it enjoys a tame behavior (i.e. splits as a finite graph of groups with finitely generated edge and vertex groups, in a way which is controlled by the underlying graph, in the sense of Theorem \ref{thm:generaldichotomy}).   \par

In this section we show that when $\Gamma$ is chordal every discrete character is subject to a complete tame-wild dichotomy; on the other hand   when $\Gamma$ is not chordal it is possible to construct characters whose kernels are not wild but at the same time are not tame in the above sense. We begin by recalling some definitions.

\begin{definition}
A graph $\Gamma$ is \textit{chordal} if for every cycle of length at least four there is a \textit{chord}, i.e. an edge connecting two non-adjacent vertices of the cycle. \end{definition}

These graphs are also known as triangulated graphs. By \cite[Theorem 1]{D87}   chordality is also equivalent to the fact that $\raag \Gamma$ is coherent (i.e. every finitely generated subgroup is finitely presented).
Dirac characterized chordal graphs in terms of their splittings, in a way that we now review.
\begin{definition}
Let $\Gamma$ be a graph and $a,b\in \vv\Gamma$. A subgraph $\Lambda \subset \Gamma$ is said to be an $ab$-\textit{separator} if it is separating and $a$ and $b$ belong to different connected components of $\Gamma \setminus \Lambda$. We say $\Lambda$ is a \textit{minimal} $ab$-\textit{separator} if it is minimal (with respect to inclusion) among $ab$-separators. We say that $\Lambda$ is a \textit{minimal vertex separator} if it is a minimal $ab$-separator for some $a,b \in \vv \Gamma$.
\end{definition}

We explicitly notice that a minimal separating subgraph (i.e. a subgraph of $\Gamma$ which is separating in the sense of Definition \ref{def:separating}, and has no proper subgraph which separates $\Gamma$) is automatically a minimal vertex separator, but the converse is not true in general, i.e. a minimal vertex separator of $\Gamma$ can contain a proper subgraph which separates $\Gamma$.
According to a classical theorem by Dirac, a graph is chordal if and only if  every minimal vertex separator is a clique (i.e. a complete graph). Notice that being chordal is an hereditary property, i.e. if $\Lambda$ is a connected full subgraph of a chordal graph $\Gamma$ then it is chordal too. In particular if $\Gamma$ is chordal and $(\Gamma_1,\Gamma_2,\Gamma_3)$ is a   splitting, then each of the subgraphs appearing in it are chordal too. Using this characterization we obtain the following result.

\begin{theorem}\label{thm:chordal}
Let $\Gamma$ be a chordal graph and let $f:\raag \Gamma \to \rr$ be discrete. Then exactly one of the following holds.
\begin{enumerate}
    \item (wild) $\ker(f)$ surjects onto $\ff_\infty$;
    
    \item (tame) $\ker(f)$ splits as a finite graph of groups $\widehat {\mathcal G}$  with finitely generated free abelian edge and vertex groups. 
\end{enumerate}
\end{theorem}
\proof
Suppose $\ker(f)$ does not surject to $\ff_\infty$; then by Corollary \ref{cor:wild} we know that $\ker(f)$ is finitely generated, and that for any separating subgraph $\Lambda$ we have that $\Lambda\not \subseteq \mathcal D(f)$. If $\Gamma$ is complete, then the statement is trivial, because there is no separating subgraph and any subgroup of $\raag \Gamma$ is already a finitely generated free abelian group.\par

So assume $\Gamma$ is not complete, in which case it admits a splitting $(\Gamma_1,\Gamma_2,\Gamma_3)$. As usual let us denote by $f_k$ the restriction of $f$ to $\Gamma_k$ for $k=1,2,3$. Up to replacing $\Gamma_3$ with a proper subgraph which still separates $\Gamma$, we can assume that $\Gamma_3$ is a minimal separator; since $\Gamma$ is chordal, $\Gamma_3$ is a clique.
As observed above, it is not possible that $f=0$ on the entire $\Gamma_3$; in particular $\mathcal L(f_3)$ is connected and dominating in $\Gamma_3$, and therefore $\ker(f_3)$ is finitely generated. By Lemma \ref{lem:BNSoversplit} we get that $\ker(f_k)$ is finitely generated for $k=1,2$ as well. So the graph of groups $\mathcal G$ coming from Theorem \ref{thm:generaldichotomy} is finite and has finitely generated vertex and edge groups. Moreover $\Gamma_1$ and $\Gamma_2$ are still chordal graphs, so we can iterate the same argument on the vertex groups of $\mathcal G$. Notice that at each step the edge groups are finitely generated free abelian groups, because we use splittings along cliques; moreover the procedure halts when the pieces of the splitting are complete graphs, and the terminal graph of groups $\widehat {\mathcal G}$ satisfies the requirements.
\endproof

\begin{remark}\label{rmk:vgbs}
Groups that can be realized as the fundamental group of a finite graph of groups    with finitely generated free abelian vertex groups have   been studied by Beeker under the name of generalized Baumslag-Solitar groups of variable rank (vGBS groups). Beeker has studied JSJ-decompositions of such groups in \cite{B13,B14}, and the multiple conjugacy problem for them in \cite{B15}. Notice that it follows from the above proof that  the vertex groups of  $\widehat {\mathcal G}$  correspond to maximal cliques of $\Gamma$, while edge groups correspond to the minimal vertex separators of $\Gamma$; in particular their rank can be easily computed. See Section \ref{sec:block_graphs} below for a detailed discussion of this point of view in the case of block graphs.
\end{remark}

For completeness we show that when the graph is not chordal it is possible to construct characters that do not obey the tame-wild dichotomy. 

\begin{proposition}\label{prop:nonchordal}
Let $\Gamma$ be connected but not chordal. Then there exists a discrete character $f:\raag \Gamma \to \rr$  such that $\ker(f) $ is neither wild nor tame. More precisely:
\begin{enumerate}
    \item (not wild): $\ker(f) $ does not surject to $\ff_\infty$,
    \item (not tame): there exists a splitting of $\Gamma$ such that the associated graphical splitting of $\ker(f) $ does not have finitely generated edge groups.
\end{enumerate}
\end{proposition}

\proof  
Since $\Gamma$ is not chordal, there exists a splitting $(\Gamma_1,\Gamma_2,\Gamma_3)$ in which $\Gamma_3$ is a non-complete minimal separating subgraph. Since $\Gamma_3$ is minimal, no proper subgraph $\Lambda\subset \Gamma_3$ can separate $\Gamma$. We distinguish two cases, depending on the connectedness of $\Gamma_3$.

If $\Gamma_3$ is disconnected, then  just set $f=0$ on one component $\Lambda$ of $\Gamma_3$, and $f=1$ on any other vertex of $\Gamma$. Then by construction $\Lambda = \mathcal D(f)$ is a proper subgraph of $\Gamma_3$. Moreover $\mathcal L(f_3)$ is not dominating in $\Gamma_3$, hence $\ker(f_3)$ is not finitely generated and the graph of groups $\mathcal G$ from Theorem \ref{thm:generaldichotomy} does not have finitely generated edge groups. Since $\Gamma_3$ was chosen to be  a minimal separating subgraph, no subgraph of $\Lambda = \mathcal D(f)$ can separate $\Gamma$. Then by Corollary \ref{cor:wild} $\ker(f)$ is not wild.

So let $\Gamma_3$ be connected; since it is not complete, it admits a separating subgraph $\Lambda$.
Then set $f=0$ on each vertex of $\Lambda$ and $f=1$ on any other vertex of $\Gamma$. Once again,   $\Lambda=\mathcal D(f)$ can not contain any  subgraph that separates $\Gamma$ by minimality of $\Lambda$, so $\ker(f)$ is not wild by Corollary \ref{cor:wild}. The same Corollary can be applied to $\Gamma_3$ to argue that since $\Lambda$ separates $\Gamma_3$, we have that $\ker(f_3)$ is not finitely generated, which implies the statement as above.
\endproof

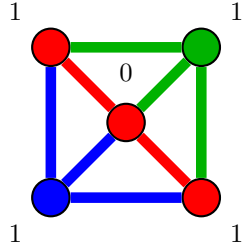
\begin{figure}[h]
    \centering
    \begin{tikzpicture}[thick]
        \tikzstyle{every node}=[circle, draw, fill=black!50,
                            inner sep=5pt, minimum width=2pt]
        
        
            \node (n0)[fill=blue, label=-135:{$1$}]  at (0,0) {};
            \node (n1)[fill=red, label=135:{$1$}]  at (0,2) {};
            \node (n2)[fill=red, label=above:{$0$}]  at (1,1) {};
            \node (n3)[fill=red, label=-45:{$1$}]  at (2,0) {};
            \node (n4)[fill=green!70!black, label=45:{$1$}]  at (2,2) {};
            
        \path[every node/.style={font=\sffamily},line width=4pt]
        
            (n0) edge[blue] node {} (n1)
            (n0) edge[blue] node {} (n3)
            (n1) edge[green!70!black] node {} (n4)
            (n2) edge[blue] node {} (n0)
            (n2) edge[red] node {} (n1)
            (n2) edge[red] node {} (n3)
            (n2) edge[green!70!black] node {} (n4)
            (n3) edge[green!70!black] node {} (n4);
        
    \end{tikzpicture}
    \caption{A non-chordal well-connected graph with a character as in Proposition \ref{prop:nonchordal}}\label{fig:nonchordalwellconnected}
\end{figure}
 
Notice that the statement of the above Proposition \ref{prop:nonchordal} is somewhat obvious for graphs admitting a non-connected splittings; for instance if $\Gamma$ is a  cycle of at least four vertices, then the Bestvina-Brady group is finitely generated but not finitely presented, and the only graphical splitting available is such that the edge group is not finitely generated. On the other hand Proposition \ref{prop:nonchordal} also applies to graphs for which every separating subgraph is connected (but not necessarily complete), such as well-connected graphs (see \cite{DLVM86}); these include for instance  the cone or the suspension over any graph (see Figure \ref{fig:nonchordalwellconnected} for an explicit example).

\begin{remark}
We conclude this section by explicitly observing that when we are in the tame case we obtain a decomposition of the group as a finite graph of groups in a way that can be computed from the underlying graph $\Gamma$, which is what we would call a graphical splitting.
It is not clear whether this is the only way in which an Artin kernel can split as a finite graph of groups with finitely generated vertex and edge groups. Obtaining such a result is outside the scope of the present work.
It should be noted however that Chang has recently shown in \cite{C20} that the Bestvina-Brady group on $\Gamma$ splits  over abelian subgroups if and only if $\Gamma$ is disconnected, complete, or has a separating clique.
\end{remark}

\subsection{Block graphs}\label{sec:block_graphs}

In this section we take a more computational point of view, and obtain explicit decompositions and rank formulas for tame Artin kernels over graphs with cut vertices, in the spirit of Remark \ref{rmk:vgbs}. For simplicity, we will directly work with integral characters $f:\raag \Gamma \to \zz$.

\begin{definition}
A graph $\Gamma$ is \textit{biconnected} if it is connected and it does not have any cut vertices. The \textit{blocks} of $\Gamma$ are its biconnected components (i.e. its maximal biconnected full subgraphs); we denote by $\block{\Gamma}$ the set of blocks of $\Gamma$. 
\end{definition}

In the following we write $\ind_\Gamma=[\zz:f(\raag{\Gamma})]$; moreover for each block $B\in \block{\Gamma}$ we denote by $f_B$ the restriction of $f$ to the subgroup $\raag B$, and we define  $\ind_{\Gamma,B}=[f(\raag{\Gamma}):f(\raag{B})]$. 
The block degree $\bldeg_\Gamma(v)$ of a vertex $v \in \vv \Gamma$ is defined to be the number of blocks of $\Gamma$ meeting at $v$. Notice that if $v$ is not a cut vertex, then $\bldeg_\Gamma(v)=1$. We find it convenient to introduce the following number associated to  $\Gamma$ and to  $f:\vv \Gamma \to \zz$:

$$m(\Gamma,f)= 1-\sum_{B\in \block{\Gamma}} \ind_{\Gamma,B}+\sum_{v\in \vv \Gamma}\left( \bldeg_\Gamma(v)-1 \right)\dfrac{|f(v)|}{\ind_\Gamma} $$

Notice that it follows from the above general discussion that if $f$ vanishes on a cut vertex, then we get a surjection to $\ff_\infty$.  When this is not the case, we can get the following result.

\begin{proposition}\label{prop:split_blocks}
Let $\Gamma$ be a connected graph, and $f:A_\Gamma\to \mathbb Z$ be non-zero on cut vertices. Then 

$$ \ker(f)\cong \mathbb F_{m(\Gamma,f)} \ast \left( \underset{B\in \block{\Gamma}}{\bigast} \left(\ker \left( f_B \right) \right)^{*\ind_{\Gamma,B}} \right) $$

\end{proposition}
\proof 

We argue by induction on the number of blocks of $\Gamma$. The theorem is trivial in the case in which $\Gamma$ is biconnected, as $\Gamma$ consists of a single block and $m(\Gamma,f)=0$. So let us assume there are at least two blocks (equivalently, at least a cut vertex).\par

Let us pick a block $B_0\in \block{\Gamma}$ which contains only one of the cut vertices of $\Gamma$ (this exists, because $\Gamma$ is now assumed to be non-biconnected). Then let  $\Lambda=\Gamma / B_0$ be the quotient graph, and let $v_0$ be the unique cut vertex of $\Gamma$ separating $B_0$ from $\Lambda$. Notice $(\Lambda,B_0,v_0)$ is a connected splitting of $\Gamma$
Since $f(v_0)\neq 0$, all the indices $\ind_{\Gamma,\Lambda},\ind_{\Gamma,v_0}$ and $\ind_{\Gamma,B_0}$ are non-zero, hence by Theorem \ref{thm:generaldichotomy} $\ker(f)$ decomposes as a finite graph of groups $\mathcal G$; moreover $\mathcal G$ has $\ind_{\Gamma,\Lambda}$ vertices with vertex groups isomorphic to $\ker \left( f_\Lambda  \right)$, $\ind_{\Gamma,B_0}$ vertices with vertex groups isomorphic to $\ker \left( f_{B_0} \right)$, and $\ind_{\Gamma,v_0}=\frac{|f(v_0)|}{\ind_\Gamma}$ edges with trivial edge group. Since all edge groups are trivial, in particular  we get

$$ \ker(f) \cong  \ff_{b_1} \ast  \left(\ker \left( f_\Lambda \right) \right)^{*\ind_{\Gamma,\Lambda}}  \ast  \left(\ker \left( f_{B_0} \right) \right)^{*\ind_{\Gamma,B_0}} $$

where $b_1=1-\left( \ind_{\Gamma,\Lambda}+\ind_{\Gamma,B_0}\right) + \dfrac{|f(v_0)|}{\ind_\Gamma}$ is the first Betti number of the underlying topological  graph of $\mathcal G$.
Now notice that $\block{\Gamma}=\{B_0\}\cup \block\Lambda$, and each cut vertex of $\Lambda$ is also a cut vertex in $\Gamma$ (even though  $v_0$ might be a cut vertex in $\Lambda$ or not). So by induction we get

$$ \ker \left( f_\Lambda \right)\cong \ff_{m(\Lambda,f)} \ast \left( \underset{B\in \block \Lambda}{\bigast} \left(\ker \left( f_B \right) \right)^{*\ind_{\Lambda,B}} \right) $$

Therefore we obtain

$$ \ker(f) \cong  \ff_{b_1} \ast  \left(\ff_{m(\Lambda,f)} \ast \left( \underset{B\in \block \Lambda}{\bigast} \left(\ker \left( f_B \right) \right)^{*\ind_{\Lambda,B}}  \right) \right)^{*\ind_{\Gamma,\Lambda}}  \ast  \left(\ker \left( f_{ B_0} \right) \right)^{*\ind_{\Gamma,B_0}} \cong$$

$$\cong  \ff_{b_1+\ind_{\Gamma,\Lambda}m(\Lambda,f)} \ast  \left( \underset{B\in \block \Gamma}{\bigast} \left(\ker \left( f_B \right) \right)^{*\ind_{\Gamma,B}} \right) $$

where we used that $\ind_{\Gamma,\Lambda}\ind_{\Lambda,B}=\ind_{\Gamma,B}$ for all $B\in \block \Gamma$. We are therefore left with checking that the rank of the free group is the right one, i.e. with showing that $b_1+\ind_{\Gamma,\Lambda}m(\Lambda,f)=m(\Gamma,f)$. This can be checked by a direct computation, using the relation $\ind_{\Gamma,\Lambda}\ind_{\Lambda,B}=\ind_{\Gamma,B}$  once again, and also that $\ind_\Lambda=\ind_\Gamma \ind_{\Gamma,\Lambda}$, together with the fact that $\bldeg_\Gamma(v)=\bldeg_\Lambda(v)$ for all vertices $v\neq v_0$, while $\bldeg_\Lambda(v_0)=\bldeg_\Gamma(v_0)-1$; the computation is as follows:

$$b_1+\ind_{\Gamma,\Lambda}m(\Lambda,f)=$$

$$=1-\left( \ind_{\Gamma,\Lambda}+\ind_{\Gamma,B_0}\right) + \dfrac{|f(v_0)|}{\ind_\Gamma}+ \ind_{\Gamma,\Lambda}-\sum_{B\in \block{\Lambda}} \ind_{\Gamma,\Lambda}\ind_{\Lambda,B}+\sum_{v\in \vv \Lambda}\left( \bldeg_\Lambda(v)-1 \right)\dfrac{|f(v)|\ind_{\Gamma,\Lambda}}{\ind_\Lambda}=$$

$$= 1- \sum_{B\in \block{\Gamma}} \ind_{\Gamma,B} +\sum_{v\in \vv \Gamma}\left( \bldeg_\Gamma(v)-1 \right)\dfrac{|f(v)|}{\ind_\Gamma}=m(\Gamma,f)$$

\endproof

\begin{remark}
The fact that $\Gamma$ can be decomposed into its blocks by cutting along its cut vertices, allows for a hierarchical decomposition of $\raag \Gamma$. In the terminology of \cite{CL16} this would be called a good $\zz$-hierarchy, with leaf groups corresponding to the blocks of $\Gamma$ (i.e. the leaf groups are the RAAGs $\raag B,B \in \block{\Gamma}$). Notice that a RAAG is slender if and only if its defining graph is complete: indeed it contains $\ff_2$ otherwise. From this point of view, our Theorem \ref{prop:split_blocks} can be seen as a generalization of \cite[Theorem 4.4 (2)]{CL16} to this setting. In order for the aforementioned hierarchy to have slender leaf groups, one needs to work with graphs such that each block is a complete graph, which is the case on which we focus next.
\end{remark}

\begin{definition}
A \textit{block graph} is a graph whose blocks are complete subgraphs.
\end{definition}

Notice that in this case blocks are automatically the maximal cliques of $\Gamma$. We obtain the following statement, involving an explicit rank formula, as a direct consequence of Theorem \ref{prop:split_blocks}.

\begin{theorem}\label{thm:rankformula}
Let $\Gamma$ be a connected block graph, and $f:A_\Gamma\to \mathbb Z$ be non-zero on cut vertices. Then
$$ \ker(f)\cong \mathbb F_{m(\Gamma,f)} \ast \left( \underset{B\in \block{\Gamma}}{\bigast} \left(\zz^{|\vv B|-1} \right)^{*\ind_{\Gamma,B}} \right) $$
In particular
$$\rk( \ker(f))= 1+  \sum_{B\in \block \Gamma} \ind_{\Gamma,B}(|\vv B|-2)  +\sum_{v\in \vv \Gamma}\left( \bldeg_\Gamma(v)-1 \right)\dfrac{|f(v)|}{\ind_\Gamma} $$
\end{theorem}
\proof
It is enough to notice that when $B$ is a clique with $|\vv B|$ vertices we have that $\raag B\cong \zz^{|\vv B|}$ and therefore $\ker(f_B)\cong \zz^{|\vv B|-1}$. Computing the rank results in:
$$\rk(\ker(f))= 1-\sum_{B\in \block{\Gamma}} \ind_{\Gamma,B}+\sum_{v\in \vv \Gamma}\left( \bldeg_\Gamma(v)-1 \right)\dfrac{|f(v)|}{\ind_\Gamma} + \sum_{B\in \block \Gamma} \ind_{\Gamma,B}(|\vv B|-1) $$
\endproof

\begin{remark}\label{rmk:CL} 
When $\Gamma$ is a tree, the group $\raag \Gamma$ is free-by-cyclic; indeed it follows from \cite{D87} that the Bestvina-Brady group associated to a tree is a free group generated by the edges of $\Gamma$ and the RAAG $\raag \Gamma$ belongs to the class of groups studied by Cashen and Levitt in \cite{CL16}. 
When $\Gamma$ is a more general block graph,  the RAAG $\raag \Gamma$ still admits a natural good $\zz$-hierarchy with free abelian leaf groups  coming from the blocks. The above statement can be seen as a generalization of \cite[Theorem 6.1]{CL16} to more general block graphs: the elements $t_i$ which appear in the  main statements of \cite{CL16} can be identified to be the cut vertices of the underlying graph; their multiplicities correspond to their block degrees.
\end{remark}

We will now discuss some applications of the above rank formula to the study of fibrations of a RAAG $\raag \Gamma$ associated to a block graph $\Gamma$, i.e. the collection of tame Artin kernels of $\raag \Gamma$. In the spirit of Remark \ref{rmk:CL}, the Corollaries \ref{cor:minimal_fibers}, \ref{cor:unbounded_fibers} and \ref{cor:bounded_div} should be compared  to \cite[Corollary 6.3, Corollary 6.4, Example 5.14]{CL16}.

\begin{corollary}[Minimal rank]\label{cor:minimal_fibers}
If $\Gamma$ is a connected block graph, then
\begin{enumerate}
    \item the minimal rank of the kernel of a non-zero homomorphism   $f:\raag \Gamma \to \zz$ is $$\mu_{\Gamma}=1+  \sum_{B\in \block \Gamma} (|\vv B|-2)  +\sum_{v\in \vv \Gamma}\left( \bldeg_\Gamma(v)-1 \right) $$
    \item $\rk(\ker(f))=\mu_\Gamma$ if and only if $f(v)=\ind_\Gamma$ for each cut vertex.
\end{enumerate}
\end{corollary}
\proof  
Let $f:\raag \Gamma \to \zz$ be a non-zero homomorphism.
Recall that a necessary condition for $f$ to have a finitely generated kernel is that it does not send any cut vertex to $0$. Hence $|f(v)|\geq 1$ for each cut vertex $v$; in particular $\ind_\Gamma$ is finite. Moreover notice that for every blocks $B$ and cut vertex $v\in \vv B$ we have $\ind_\Gamma \leq \ind_B \leq |f(v)|$, i.e. $1 \leq \ind_{\Gamma,B} \leq \frac{|f(v)|}{\ind_\Gamma}$.
Then the rank formula implies that $\mu_\Gamma$ is a lower bound for the rank of $\ker(f)$, and that this lower bound is realized by characters that send all cut vertices to $\ind_\Gamma$.
\endproof

\begin{corollary}[Unbounded rank]\label{cor:unbounded_fibers}
If $\Gamma$ is a connected block graph with at least one cut vertex, then $A_\Gamma$ admits integral characters with kernel of arbitrarily large rank.
\end{corollary}
\proof
For each $n\in \nn$ let $f_n$ be the character that sends all cut vertices to $n$ and the other vertices to $1$. Then we have $\ind_{\Gamma,B}=\ind_\Gamma=1$, hence
$$\rk(\ker(f_n))=1+  \sum_{B\in \block \Gamma} (|\vv B|-2)  +n\sum_{v\in \vv \Gamma}\left( \bldeg_\Gamma(v)-1 \right) \to \infty$$
\endproof

For the next statement recall from \S \ref{sec:characters} that $\chars{\raag {\Gamma}}$ denotes the character sphere of $\raag \Gamma$. This result can be interpreted as an example of exotic divergence in the BNS invariant of $\raag \Gamma$, which shows the non-properness of the function that records the rank of the kernel.

\begin{corollary}[Bounded divergence]\label{cor:bounded_div}
Let $\Gamma$ be a connected block graph with at least one cut vertex, and let $f_n,f_\infty$ be integral characters such that $[f_n]\to [f_\infty]$ in $\chars{\raag {\Gamma}}$. Then the following are equivalent:
\begin{enumerate}
    \item $\exists C >0$ such that $\forall n$ we have $\rk(\ker(f_n))\leq C$, but $\ker(f_\infty)$ is not finitely generated;
    \item $\exists C' >0$ such that for all $n$ and all cut vertices $v\in \vv \Gamma$ we have $|f_n(v)|\leq C'$, but $\exists w \in \vv \Gamma$ such that $|f_n(w)|\to \infty$.
\end{enumerate}
Moreover in this situation $f_\infty(v)=0$ for all cut vertices $v\in \vv \Gamma$.
\end{corollary}
\proof
Suppose (2) holds, and consider the character $g_n=\dfrac{f_n}{|f_n(w)|}$. Then $\ker(f_n)=\ker(g_n)$ and $[g_n]=[f_n]\to [f_\infty]$. Since we have $\ind_{\Gamma,B}\leq |f_n(v)|\leq C'$ for each cut vertex $v\in \vv B$, the rank formula shows that there exists some $C>0$ such that $\rk(\ker(f_n))\leq C $ uniformly on $n$. However  for each cut vertex $v$ we have that $g_n(v)\to 0$, so that $f_\infty(v)=0$ and $\ker(f_\infty)$ can not be finitely generated.\par

Suppose (1) holds. Then by the rank formula $|f_n(v)|\leq C$ for all cut vertices. If $f_n$ was uniformly bounded on the entire $\vv \Gamma$, then $f_n$ would converge  to a limit character with finitely generated kernel; hence necessarily $f_n$ must diverge on at least one vertex, which can not be a cut vertex. As a result also in this case, up to rescaling, $f_n(v)\to 0$ for all cut vertices, so that $f_\infty(v)=0$.

\endproof

For  a concrete example which displays this behavior see the following.

\begin{example}\label{ex:bounded_div}
Let $\Gamma$ be a connected block graph with at least one cut vertex. To simplify the rank computation, also assume that each block contains at least one vertex which is not a cut vertex (the following computation can be adapted to remove this extra assumption).
Pick a number $q \in \zz,q\neq 0$ and an infinite sequence of numbers $p_n\in \zz$ diverging to infinity and each coprime with $q$. 
Consider the character $f_n$ defined by sending all cut vertices of $\Gamma$ to $q$ and all other vertices to $p_n$. We have $\ind_\Gamma=1$ and, since every block contains a vertex which is not a cut vertex, we also have $\ind_{\Gamma,B}=1$  for each block. As a result
$$\rk(\ker(f_n))= 1+ \sum_{B\in\block \Gamma}(|\vv B|-2)+q\sum_{v\in \vv \Gamma}(\bldeg_\Gamma(v)-1)$$
does not depend on $n$ (indeed by Proposition \ref{prop:split_blocks} all these groups are actually isomorphic).
Also consider the character $g_n=\dfrac{f_n}{p_n}$, i.e. the one obtained by sending all cut vertices to $\dfrac{q}{p_n}$ and all other vertices to $1$. We have $[f_n]=[g_n]\in \chars{\raag \Gamma}$. However $g_n\to g_\infty$, where $g_\infty$ is a character defined by sending all cut vertices to $0$ and all other vertices to $1$; by Theorem  \ref{thm:living_subgraph} $\ker(g_\infty)$ is not finitely generated.
\end{example}

\begin{remark}
It is natural to ask for a deeper understanding of this phenomenon, i.e. what makes these subgroups, which are in some explicit sense uniformly tame, abruptly converge to a wild subgroup. In terms of the graphical splittings studied in this paper it can be noticed that these subgroups are all isomorphic as vGBS groups, i.e. they can be realized as the fundamental group of a fixed graph of free abelian groups $\widehat{\mathcal G}$, determined by the structure of $\Gamma$ itself. In other words no intrinsic algebraic invariant of these subgroups seem to be useful in understanding this phenomenon. 
From a more geometric point of view, one could try to understand the extrinsic properties of these as subgroups of $\raag \Gamma$, e.g. if  these tame subgroups are more and more distorted in $\raag \Gamma$ as they approach the wild limit.  In this direction it should be noticed that the work of Tran  shows that any finitely generated Artin kernel has at worst quadratic distortion (see \cite[Corollary 1.5]{T16}), and at worst quadratic relative divergence (see \cite[Proposition 4.3]{T17}).
\end{remark}

\printbibliography

\end{document}